\documentclass[12pt]{article}
\usepackage{latexsym}
\begin{document}
\thispagestyle{empty}

\vskip 20pt
\begin{center}
{\bf PERMUTATIONS RESTRICTED BY TWO DISTINCT\\ PATTERNS OF LENGTH
THREE}
\vskip 15pt
{\bf Aaron Robertson}\\
{\it Department of Mathematics,}
{\it Colgate University,
Hamilton, NY 13346}\\
{\tt aaron@math.colgate.edu}
\end{center}
\vskip 30pt
\begin{abstract}{\footnotesize \noindent
Define $S_n(R;T)$ to be the number of permutations on $n$ letters
which avoid all patterns in the set $R$ and contain each
pattern in the multiset $T$ exactly once.  In this paper we enumerate
$S_n(\emptyset;\{\alpha,\beta\})$ for
all
$\alpha \neq \beta \in S_3$.  
}
\end{abstract}
\vskip 20pt
{\bf 1. Introduction}

Let $\pi \in S_n$ be a permutation of $[n]=\{1,2,\dots,n\}$ written as a
word.
Let $\alpha \in S_k$, $k \leq n$.
We say that $\pi$ {\it contains the pattern $\alpha$} if there exist
indices $i_1,i_2,\dots,i_k$ such that $\pi_{i_1} \pi_{i_2} \dots
\pi_{i_k}$ is equivalent to $\alpha$, where we define equivalence as
follows. Define $\overline{\pi}_{i_j}=
|\{m:\pi_{i_m} \leq \pi_{i_j}, m=1,2,\dots,k\}|$.  If
$\alpha = \overline{\pi}_{i_1} \overline{\pi}_{i_2} \dots
\overline{\pi}_{i_k}$ then
we say that $\alpha$ and $\pi_{i_1} \pi_{i_2} \dots \pi_{i_k}$ are
equivalent.  For example, if $\tau=124635$ then $\tau$ contains
the pattern $213$ by noting that $\tau_3 \tau_5 \tau_6 = 435$
is equivalent to $213$.  We say that $\pi$ {\it avoids the pattern
$\alpha$} if $\pi$ does not contain the pattern $\alpha$.
In our above example, $\tau$ avoids the pattern $321$.

Let $\alpha \neq \beta$ be patterns of length three.  In this
article we enumerate the number of permutations which contain
$\alpha$ exactly once and avoid $\beta$ as well as those
permutations which contain each of $\alpha$ and $\beta$
exactly once.

{\bf 2. Some History}

The investigation of permutations which avoid a pattern of
length three started well over a hundred years ago as
exhibited in [C] and references therein.
Knuth ([Kn]) investigated permutations which avoid any
single
pattern of length $3$ and showed that, regardless of the pattern,
such permutations are enumerated by the Catalan numbers.
Bijective results are given in [Ri], [Krt], [SS], and [W1].  To describe
the enumeration results more succinctly we introduce the following
notation. Let $S_n(R)$ be the set of permutations on $[n]$
which avoid all patterns in the set $R$, where we omit
the set notation if $|R|=1$, and let
$s_n(R) = |S_n(R)|$.  Knuth's result
can then be stated as $s_n(\alpha) = \frac{1}{n+1} {2n \choose n}$ for
all $\alpha \in S_3$.

Following Knuth's result, two natural progressions
were made:  the investigation of $S_n(R)$ for $R \subseteq S_3$ and the
investigation of
$S_n(\beta)$ for $\beta \in S_4$.
With respect to the former investigation,
Simion and Schmidt ([SS]) gave a complete study of
$s_n(R)$ for all $R \subseteq S_3$.
With respect to the latter investigation, in two beautiful
papers, Gessel ([Ge]) found
$s_n(1234)$ and B\'ona ([B1]) found
$s_n(1342)$.
Further results on $S_n(\alpha)$ for $\alpha \in S_4$
are given by West in
[W1] and [W2] and by Stankova in [S].
The exact enumeration of
$1324$-avoiding permutations is still an open question, with
the only result being a lower bound given by
B\'ona in [B2].

Several logical extensions followed:  the investigation of
$S_n(R)$ for $R \subseteq S_4$, the investigation of
$S_n(S \cup T)$ for $S \subseteq S_3$ and $T \subseteq S_4$,
and the investigation of $S_n(R)$ for $R \subseteq S_j$, $j>4$.
Guibert, in [Gu], showed that for certain $R \subseteq S_4$
with two elements, the corresponding $s_n(R)$ are given by
Schr\"oder
numbers.  In [B3] and [Kr], B\'ona and Kremer, respectively, gave
further
extensions for
$R \subseteq S_4$ with two elements.
Mansour ([M])
completely enumerated $S_n(R \cup \{\alpha\})$ for $R \subseteq S_3$
and $\alpha \in S_4$.
Results for permutations avoiding patterns of length greater than four
can
be found in  [BLPP1], [BLPP2], [CW], and [Kr].

A natural generalization of pattern-avoiding permutations
is pattern-containing permutations.  To aid in the discussion
of pattern-containing permutations we introduce the following
notation.
Let $S_n(R;T)$ be the set of permutations on $[n]$
which avoid all patterns in the set $R$
and contain each
pattern in the multiset $T$ exactly once, where we again omit
the set notation for singleton sets, and let
$s_n(R;T) = |S_n(R;T)|$.

Recently, there has been much research focused on
$S_n(R;T)$ for various sets $R$ and multisets
$T$.  Below, we give some
results in this direction.
First, in [N], Noonan proved that
$s_n(\emptyset;123)=\frac{3}{n} {2n \choose n+3}$, a remarkably
elegant formula.  B\'ona, in [B4], then showed that
$s_n(\emptyset;132) = {2n-3 \choose n-3}$, an even simpler
formula, proving a conjecture presented in [NZ].  These two results give
$s_n(\emptyset;\alpha)$ for all $\alpha \in S_3$, by applying
the following two bijections (given in [SS]).

{\bf Reversal}:  Define $r: S_n \rightarrow S_n$ by
$r(\pi_1 \pi_2 \dots \pi_n) = \pi_n \pi_{n-1}
\dots \pi_1$.

{\bf Complementation}:  Define $c:  S_n \rightarrow S_n$ by
$c(\pi_1 \pi_2 \dots \pi_n) = (n-\pi_1+1)
(n-\pi_2+1) \dots (n-\pi_n+1)$

We will also have need of a third bijection (given in [SS])
which is defined as follows.

{\bf Inverse}:  Define $i: S_n \rightarrow S_n$
as the group theoretic inverse.

It is easy to see that if $\pi$ contains exactly $s \geq 0$
occurences of the pattern
$\alpha$, then $r(\pi)$ (resp. $c(\pi)$, $i(\pi)$) contains
exactly
$s$ occurences of the pattern
$r(\alpha)$ (resp. $c(\pi)$, $i(\pi)$).  By applying
$r$, $c$, and $r \circ c$ we see that
$s_n(\emptyset;123)=s_n(\emptyset;321)$ and
$s_n(\emptyset;132)=s_n(\emptyset;231)
=s_n(\emptyset;312)=s_n(\emptyset;213)$.

In [B4], B\'ona also gave the generating function for
$\{s_n(\emptyset;\{132,132\})\}_n$.  In [R], the formulas for
$s_n(132;123), s_n(123;132)$, and $s_n(\emptyset;\{123,132\})$
are given.  These results were extended in [RWZ] to give the
generating function for $\{s_n(132;\{123^r\})\}_{r,n \geq 0}$
in the form of a continued fraction.  Mansour and Vainshtein
([MV1]) generalized this result to give the
generating function for
$\{s_n(132;\{(123 \dots k)^r\})\}_{r,n}$ for a
given $k$
and showed the relation of such permutations to
Chebyshev polynomials of the second kind.
In [CW] other similar permutations were first
shown to be related to the Chebyshev polynomials of
the second kind.
Independently, Jani and Rieper ([JR]) also extended the
result in [RWZ] to find the generating function given in
[MV1] using the theory of ordered trees.  Shortly thereafter,
Krattenthaler, in [Krt], used Dyck path bijections to reprove elegantly
the results in [MV1] and [JR], extend results given in
[CW], give a precise asymptotic formula for $s_n(132,\{(123 \dots
k)^r\})$,
and show that
$s_n(132,\{(123 \dots k)^r\}) \asymp s_n(123,\{((k-1)(k-2)\dots 1
k)^r\})$

\section*{\normalsize 3. Preliminaires}

In this section we give some definitions
and state a known result (without proof) upon which we will need to
draw.

In order to discuss our analysis we have need
of the following two definitions.  The first definition
has become a standard definition, while the second
definition is new.

{\bf Definition} (Wilf class)  Let $S_1$ and $S_2$ be
two sets.  If
$s_n(S_1)=s_n(S_2)$ then we say that
$S_1$ and $S_2$ are in the same {\it Wilf class}, or
{\it Wilf equivalent}.

{\it Example.}  There is only one Wilf class for
permutations avoiding a single pattern of
length $3$ since
$s_n(\alpha) = \frac{1}{n+1} {2n \choose n}$ for
any $\alpha \in S_3$.

{\bf Definition} (almost-Wilf class\footnote{
As an aside, Herb Wilf has told the author that
he is not fond of the monicker Wilf class,
however, in honor of Herb (and due to the lack
of a better name), we extend what has become the
standardized definition of pattern-avoiding
permutation classes.})  Let $S_1$ and
$S_2$ be two sets and let
$T_1$ and $T_2$ be two multisets.  If
$s_n(S_1;T_1) = s_n(S_2;T_2)$ then we say
that $(S_1;T_1)$ and $(S_2;T_2)$ are in the same
{\it almost-Wilf class}, or
{\it almost-Wilf equivalent}.

{\bf Theorem 2.1}  (Simion and Schmidt, [SS])

{\it
\hskip 10pt
(1). For $\{\alpha,\beta\} \in
\{ \{123,132\}, \{123,213\}, \{132,213\},
\{132,231\}, \{132,312\}, \{213,231\},$
\vskip -10pt
\hskip 34pt
$\{213,312\},\{231,312\}, \{231,321\}, \{312,321\}\}$ we have
$s_n(\{\alpha,\beta\}) =2^{n-1}$ for $n \geq$
\vskip -10pt
\hskip 34pt
 $2$ and $s_1(\{\alpha,\beta\})=1$;

\hskip 10pt
(2). For $\{\alpha,\beta\} \in
\{\{123,231\}, \{123,312\}, \{132,312\}, \{213,321\}\}$
we have
$s_n(\{\alpha,\beta\})=$
\vskip -10pt
\hskip 34pt
 ${n \choose 2}+1$;

\hskip 10pt
(3). $s_n(\{123,321\})=0$ for $n \geq 5$.}

\section*{\normalsize 4. On $s_n(\alpha;\beta)$}

As seen in Section 2 we know
$s_n(\alpha;\beta)$ for
$(\alpha,\beta)\in \{(123,132),(132,123)\}$.  Using
the reversal and complementation bijections
presented in Section 2 we see that the following
is true.

{\bf Theorem 4.1} {\it For $(\alpha,\beta) \in
\{(123,132), (123,213), (132,123), (213,123),
(231,321), (312,$
\vskip -10pt
\hskip 79pt
$321), (321,231), (321,312) \}$ we
have $s_n(\alpha;\beta)=(n-2)2^{n-3}$ for $n \geq 3$.}

To complete the enumeration $s_n(\alpha;\beta)$ for
all $\alpha \neq \beta \in S_3$
we must consider the following classes, which can
be obtained through application of the 
reversal, complementation, and
inverse bijections.

(1).  $\{(123;\!321), (321;\!123)  \}$

(2).  $\{(123,231), (123,312), (321,132), (321,213)  \}$

(3).  $\{(132;\!213), (213;\!132), (231;\!312), (312;\!231)  \}$

(4). $\{(132;\!231),(132;\!312), (213;\!231),(213;\!312), (231;\!132),
(231;\!213), (312;\!132),(312;\!213)\}$

(5).  $\{(132,321),(213,321), (231,123), (312,123) \}$

Trivially, we have $s_n(123;321)=0$ for $n \geq 6$. 
The enumeration concerning the remaining classes
follows from results, which will be noted below,
given by Mansour and Vanshtein
in [MV2] and [MV3].

{\bf Theorem 4.2} {\it For $(\alpha,\beta) \in
\{(123,231), (123,312), (321,132), (321,213)  \}$ we
have
\vskip -10pt
\hskip 79pt
$s_n(\alpha;\beta) =2n-5$
for
$n \geq 3$.}

{\it Proof.}  This follows from Theorem 3.3 in
[MV3] with $m=2$ and $k=3$.
\hfill{$\Box$}

{\bf Theorem 4.3} {\it For $(\alpha,\beta) \in
\{(132,213), (213,132), (231,312), (312,231)  \}$ we
have
\vskip -10pt
\hskip 79pt
$s_n(\alpha;\beta) = n2^{n-5}$ for
$n \geq 4$ and $s_3(\alpha;\beta)=1$.}

{\it Proof.}  This follows from Example 3.2 in [MV2]
with $p=1, m=2$ and $k=3$.
\hfill{$\Box$}

{\bf Theorem 4.4} {\it For $(\alpha,\beta) \in
\{(132,231), (132,312), (213,231), (213,312), (231,132),
(231,$
\vskip -10pt
\hskip 79pt
$213), (312,132), (312,213) \}$ we have
$s_n(\alpha;\beta) =2^{n-3}$ for
$n \geq 3$.}

{\it Proof.}  This follows from Theorem 3.4 in [MV2]
with $m=1$ and $k=3$.
\hfill{$\Box$}

{\bf Theorem 4.5} {\it For $(\alpha,\beta) \in
\{(132,321),(213,321), (231,123), (312,123) \}$ we
have
\vskip -10pt
\hskip 79pt
$s_n(\alpha;\beta) = 2n-5$ for
$n \geq 3$.}

{\it Proof.}  This follows immediately from Theorem 3.2 in [MV2].
\hfill{$\Box$}

{\it Remark.}  Notice that {\it a priori} there were
six classes we had to consider
(by Theorems 4.1 through 4.5 and the trivial case).  (This is
one less than the seven classes to consider
before [R] showed that
$(123;132)$ and $(132;123)$ are
almost-Wilf equivalent.)  However,
the results above show that
there are in fact only five almost-Wilf classes associated
with $S_n(\alpha;\beta)$, $\alpha \neq \beta \in S_3$.
An explanation of this is given in the following subsection.

\subsection*{\normalsize 4.1 Generating $S_n(123;312)$   and
$S_n(312;123)$}

In this short subsection we investigate the nature as to why
$s_n(123;312)=s_n(312;123)$ (which are both equal to $2n-5$).  We will show that the two
sets considered here are generated by almost exactly the
same rule, and let the reader infer a bijection from this
result.
Define $\phi:S_{m-1} \rightarrow S_{m}$ by $\phi(\pi_1 \pi_2 \dots \pi_{m-1})
= (\pi_1+1) (\pi_2+1) \dots (\pi_{m-1}+1) 1$.

It is easy to see that for any $\sigma \in S_{n-1}(123;312)$ and
any $\tau \in S_{n-1}(312;123)$ that $\phi(\sigma) \in S_{n}(123;312)$ 
and $\phi(\tau) \in S_{n}(312;123)$.  Since
$S_3(123;312)=\{312\}$ and $S_3(312;123)=\{123\}$ we can use the
rules below to generate $S_n(123;312)$   and
$S_n(312;123)$.

{\tt Generating Rule for $S_{n} (123;312)$}:\\
By Theorem 4.2, it is trivial to check
that $S_{n} (123;312) = \{\phi(\pi): \pi \in
S_{n-1}(123;312)\}
\cup \{ 3 1 n (n-1) (n-2) \dots 5 4 2, (n-2) (n-3) \dots 3 2 n 1 (n-1) \}$.

{\tt Generating Rule for $S_{n} (312;123)$}:\\
By Theorem 4.5, it is trivial to check
that $S_{n} (312;123) = \{\phi(\pi): \pi \in
S_{n-1}(312;123)\}
\cup \{1 (n-1) n (n-2) (n-3) \dots 3 2, (n-2) (n-1) (n-3) (n-4) \dots 2 1 n\}$.

\section*{\normalsize 5. On $s_n(\emptyset;\{\alpha,\beta\})$}

We first note that trivially $s_n(\emptyset;\{123,321\})=0$
for $n \geq 6$.  Next,
using the bijections $r$ and $c$ we have four
classes to consider:

(1).
$\overline{\{123,231\}}=\{\{123,231\},\{123,312\},\{132,321\},\{213,321\}\}$

(2).
$\overline{\{123,132\}}=\{\{123,132\},\{123,213\},\{231,321\},\{312,321\}\}$

(3). $\overline{\{132,213\}}=\{\{132,213\},\{231,312\}\}$

(4).
$\overline{\{132,231\}}=\{\{132,231\},\{132,312\},\{213,231\},\{213,312\}\}$

Class (2) was enumerated in [R] giving the following theorem.

{\bf Theorem 5.1}  {\it For $\{\alpha,\beta\} \in
\{\{123,132\}, \{123,213\}, \{231,321\},\{312,321\} \}$ we
have
\vskip -10pt
\hskip 79pt
$s_n(\emptyset,\{\alpha,\beta\}) = (n-3)(n-4)2^{n-5}$ for
$n \geq 5$.
}

Except for the proof 
of Theorem 5.4, in the proofs below we will isolate either
the element $1$ or the element $n$ in each
permutation, $\pi$.  Denote by $\pi(1)$ the elements (in order)
to the left of the isolated element, and by $\pi(2)$
the elements (in order) to the right of the isolated element.
Hence, we have $\pi=\pi(1) \, 1 \, \pi(2)$ or $\pi
=\pi(1) \, n \, \pi(2)$.
We start with class (1).

{\bf Theorem 5.2}  {\it For $\{\alpha,\beta\} \in
\{\{123,231\}, \{123,312\}, \{132,321\},\{213,321\} \}$ we
have
\vskip -10pt
\hskip 79pt
$s_n(\emptyset,\{\alpha,\beta\}) = 2n-5$ for
$n \geq 5$ and $s_4(\emptyset,\{\alpha,\beta\}) = 2$.
}

{\it Proof.}  We will use $\{123,312\}$ for our proof.
Let $f_n=s_n(\emptyset;\{123,312\})$,
let $\pi \in S_n(\emptyset;\{123,312\})$, and let $\pi_i=n$.

We have three cases to consider:  (i) the $(312)$ pattern occurs
with $n$ as the `3' and the $(12) \in \pi(2)$, (ii) the pattern
$(312) \in \pi(1)$, and (iii) the `2' in the $(312)$ pattern
is in $\pi(2)$ while $(31) \in \pi(1)$.

We start with case (i): the $(312)$ pattern occurs
with $n$ as the `3' and the $(12) \in \pi(2)$.
Let $x$ be the `1' and $y$ be the `2' in the $(312)$ pattern.

 Write $\pi
= \pi(1) \, n \, A \, x \, B \, y \, C$, where $A, B$, and $C$
represent the portions of $\pi$ in between two distinguished
elements (either $n$ and $x$, $x$ and $y$, or $y$ and the end of $\pi$).

We will first show that $A$ is empty.  Assume otherwise and let
$a \in A$.  Then either $nay$ is another $(312)$ occurence
(if
$a<y$) or $axy$ is another $(312)$ occurence (if $a>y$).
Hence, $A$ is empty. Next, we will show that $B$ must be
empty.  Assume otherwise and let $b \in B$.  Then either
$nby$ is another $(312)$ occurence (if $b<y$) or $nxb$ is
another $(312)$ occurence (if $b>y$).  Hence, $B$ must also
be empty.  Thus we may write
$\pi=\pi(1) \, n \, x \, y\,C$.

Next, we notice that for any $c \in C$ we must have $c<x$
and for any $p \in \pi(1)$ we must have $p<y$ to avoid
another $(312)$ occurence.  Furthermore, if $C$ were to contain
a $(12)$ pattern then we would have another occurence of
$(312)$ with $n$ acting as the `3'.  Hence, the elements of
$C$ must be in decreasing order and thus our $(123)$ pattern
must start in $\pi(1)$.  Similarly, the elements of $\pi(1)$
must be in decreasing order or we would have
at least two occurences of $(123)$ with both $n$ and $y$ serving
as the `3' in the $(123)$ pattern.  Hence, there exists
$r \in \pi(1)$ with $r<x$ which produces $rxy$ as our
$(123)$ pattern.  Furthermore, all other elements in $\pi(1)$
must be larger than $x$ or else we would have another
occurence of $(123)$.  Hence,
we must have $r =\pi_{i-1}$ since
the elements in $\pi(1)$ are decreasing.  However,
if $i \neq 2$ then $\pi_1 r x$ would be another $(312)$
pattern.  Thus, $i=2$.  The last piece of information
we need is that since all elements in $C$ are less than
$x$, we must have $x=n-2$.  Thus we see that our permutations
in this case are of the form $\pi=r\,n\,(n-2) \, (n-1)\,C$
with the elements of $C$ in decreasing order.
Since we have $n-3$ choices for $r$, we have $n-3$ permutations
in this case.

Next, we look at case (ii):  the pattern $(312) \in \pi(1)$.

Let $zxy$ be the $(312)$ pattern and write
$\pi=A \, z \, B \, x \, C \, y \, D \, n \, \pi(2)$.
Notice that in this case we already have our $(123)$ pattern,
namely, $xyn$.

We first show that $A$, $B$, and $C$ are empty.  Assume otherwise
and let $a \in A$, $b \in B$, and $c \in C$.  For any $a \in A$ we see
that either $ayn$ would give another $(123)$ occurence (if
$a<y$) or that $axy$ would give another $(312)$ occurence (if
$a>y$). For $b \in B$ we see that either $byn$ would give
another
$(123)$ occurence (if $b<y$) or that $bxy$ would give another
$(312)$ occurence (if $b>y$).  For $c \in C$, either $xcn$
would be another $(123)$ occurence (if $c>y$)
or $zcy$ would be another $(312)$ occurence (if $c<y$).
Hence,
$A, B$, and $C$ must all be empty so we may write
$\pi=z \, x \, y \, D \, n \, \pi(2)$.

Next, we notice that for any element in $D$ or $\pi(2)$, that
element must be less than $x$, for otherwise we would have either
another occurence of $(312)$ with $z$ and $x$ or another $(123)$
occurence with $x$ and $y$.  This restriction
gives us $z=n-1$, $x=n-3$, and $y=n-2$.  Furthermore, the elements
in $D$ must be decreasing (to avoid another $(123)$ with $n$),
and the elements in $\pi(2)$
must be decreasing (to avoid another $(312)$ with $n$).
Even further, for all $d \in D$ and all $p \in \pi(2)$ we must
have $d>p$ or else we would have another $(312)$ occurence
with $zdp$.  Hence, the elements in both $D$ and $\pi(2)$ are
determined by the position of $n$.  Since we have $n-3$ choices
for the position of $n$, we have $n-3$ permutations in this case.

Lastly, we look at case (iii):  the `2' in the $(312)$ pattern
is in $\pi(2)$ while $(31) \in \pi(1)$.

Let $zxy$ be the $(312)$ pattern and write
$\pi=A \, z \, B \, x \, C \, n \, D \, y \, E$.

We first show that $B$ and $C$ are empty.   Assume otherwise
and let  $b \in B$ and $c \in C$.  For $b \in B$, either
$bxy$ is another occurence of $(312)$ (if $b>y$) or
$zby$ is another occurence of $(312)$ (if $b<y$).
For $c \in C$, either we get two occurences of
$(123)$ with $zcn$ and $xcn$
(if $c>z$), we get another $(312)$ occurence with
$zxc$ (if $x<c<z$), or we get another $(312)$ occurence
with $zcy$ (if $ c<x$).  Hence, we may write
$\pi=A\,z\,x\,n\,D\,y\,E$.

Next, notice that the elements in $D$
must be decreasing and the elements in $E$ must be decreasing
to avoid another occurence of $(312)$ with
$n$ serving as the `3'.  Furthermore, all elements in $D$
must be greater than $z$ and all elements in $E$ must
be less than $x$ since if either of these did not
hold we would have another $(312)$ occurence.
We then see that for all $a \in A$ we must have
$x<a<y$, otherwise if $a>y$ we would obtain another $(312)$
occurence with
$x$ and $y$, and if $a<x$ we would have two occurences of $(123)$
with $axn$ and $axy$.  We also note that
$A$ must contain exactly one element since for any $a\in A$, $azn$
produces
a $(123)$ pattern and if $A$ is empty we cannot obtain
a $(123)$ occurence.  Since $A$ is not empty we now see that
$D$ must be empty to avoid another $(123)$ occurence with $a$ and
$z$. We may now write $\pi=a\,z\,x\,n\,y\,E$, where $x<a<y$.

Since all elements in $E$ must be smaller than $x$ we see that
$x=n-4$, $y=n-2$, $z=n-1$, and $a=n-3$.  Finally, since the elements
in $E$ must be decreasing we see that we only have a single permutation
in this case (provided $n \geq 5$).

Summing over all cases we have $s_n(\emptyset,\{123,312\})=2n-5$
for $n \geq 5$.
\hfill{$\Box$}

{\it Remark}.  Notice that we have the interesting result that
$s_n(123;312) = s_n(\emptyset;\{123,312\})$ and hence
$(123;312)$ and $(\emptyset;\{123,312\})$ are almost-Wilf equivalent
(for $n \geq 5$).
This is the first nontrivial case of a ``mixed restriction"
equivalence.

We now move on to class (3) and prove the following theorem.

{\bf Theorem 5.3}  {\it For $\{\alpha,\beta\} \in
\{\{132,213\},\{231,312\}\}$ we
have $s_n(\emptyset,\{\alpha,\beta\}) = $
\vskip -10pt
\hskip 79pt
$(n^2+21n-28)2^{n-9}$ for
$n \geq 7$, $s_6(\emptyset,\{\alpha,\beta\})=17$,
$s_5(\emptyset,\{\alpha,\beta\})=6$,
\vskip -10pt
\hskip 79pt
and $s_4(\emptyset,\{\alpha,\beta\})=3$.
}

{\it Proof.}  We will use $\{231,312\}$ for our proof.
Let $f_n=s_n(\emptyset;\{231,312\})$,
let $\pi \in S_n(\emptyset;\{231,312\})$, and let $\pi_i=1$.

We have three cases to consider:  (i) the pattern $(312) \in \pi(1)$,
(ii) the pattern
$(312) \in \pi(2)$, and (iii) the $(312)$ pattern
straddles 1, i.e. the `3' is in $\pi(1)$, the `2' is
in $\pi(2)$, and 1 serves as the `1' in the pattern.

We start with case (i):  the pattern $(312) \in \pi(1)$.

Let $zxy$ be our $(312)$ pattern and write
$\pi=A\,z\,B\,x\,C\,y\,D\,1\,\pi(2)$.
Note that we already have our $(231)$ pattern
with $xy1$.

We first argue that $A, B, C,$ and $D$ must all be
empty.  Assume otherwise and let
$a \in A$, $b \in B$, $c \in C$, and $d \in D$.
We start with $c \in C$.  Clearly we must have
$c>z$ to avoid another $(312)$ occurence.  However,
this produces $zc1$ which is another $(231)$
occurence.  Hence, $C$ must be empty.  Next, we
move to $b \in B$.  We see here that
either $zby$ is another occurence of $(312)$ (if
$b<y$) or $bxy$ is another occurence of $(312)$
(if $b>y$).  Hence, $B$ must also be empty.  Now,
we look at $a \in A$.  Here, either $axy$ is another
$(312)$ occurence (if $a>y$) or
both $ay1$ is another $(231)$ occurences (if $a<y$).
Lastly, for $d \in D$, either $xd1$ would be
another occurence of $(231)$ (if $d>x$) or
$xyd$ would be another occurence of $(231)$ (if
$d<x$).  Hence, we may now write
$\pi=z\,x\,y\,1\,\pi(2)$.  Since we already have
both of the required patterns we see that
$D \in S_{n-4}(\{312,231\})$.   By Theorem 2.2
we have $2^{n-5}$ permutations in this case
for $n \geq 5$, and $1$ permutation for $n=4$.

Next we look at case (ii):  the pattern $(312) \in \pi(2)$.

Let $zxy$ be our $(312)$ pattern and write
$\pi=\pi(1)\,1\,A\,z\,B\,x\,C\,y\,D$.

We first show that $B$ must be empty.  Assume
otherwise and let $b \in B$.  Then either
$zby$ is another $(312)$ (if $b<y$) or
$bxy$ is another $(312)$ (if $b>y$).  We next note that
for any $c \in C$ we must have $c>z$ to avoid another
$(312)$ occurence.  Hence, $zcy$ is a $(231)$ pattern
for any $c \in C$.  Thus, $|C| \leq 1$.

We first consider the subcase $|C|=1$.  Let $c \in C$ so
that we have both of the required patterns in our
permutation.  Write $\pi=\pi(1)\,1\,A\,z\,x\,c\,y\,D$.
Notice that for any $p\in \pi(1)$, $a\in A$,
and $d \in D$ we must have $p<a<d$.  This holds since
we must have $p<a$ to avoid another $(312)$ occurence
with $p1a$.  We then see that for any $a \in A$
we must have $a<x$ to avoid another $(231)$ occurence
with $bza$ (if $b<z$) or another $(312)$ occurence
with $bxy$ (if $b>z$).  Lastly, we note that for any
$d \in D$ we require $d>c$ to avoid $(231)$ with
$zcd$ (if $d<z$) or another $(312)$ with $zxd$ (if
$z<d<c$).  Now since our elements in $A, B$, and $D$
are either less than $x$ or greater than $c$, we
see that $y=x+1$, $z=x+2$, and $c=x+3$.

We now notice that $\pi\,1\,A$ read as a permutation
must avoid both $(231)$ and $(312)$.  Likewise,
$D$ must avoid both $(231)$ and $(312)$.
Since the value of $x$ determines the position of $x$,
by Theorem 2.2 we have
$\sum_{x=2}^{n-4} 2^{x-1} 2^{n-x-5} = (n-3)2^{n-6}$
permutations for $n \geq 6$, one permutation for
$n=5$, and none for $n \leq 4$ in this subcase.

Next, consider the subcase $|C|=0$.  Write
$\pi=\pi(1)\,1\,A\,z\,x\,y\,D$.   We have four
subsubcases to consider:

\hskip 10pt  (a) There exists a unique $d \in D$ with $d<x$.  This gives $xyd$
as our $(231)$ pattern.

\hskip 10pt  (b) There exists a unique $a \in A$ with $x<a<y$.  This gives
$azx$ as our $(231)$ pattern.

\hskip 10pt  (c)  All elements in $\pi(1)$ and $A$ are smaller than
$x$ and our $(231)$ pattern is contained
\vskip -10pt
\hskip 29pt within $\pi(1) \, 1\, A$ while $D$ avoids both patterns.

\hskip 10pt  (d)  Our $(231)$ pattern is contained within $D$
while
$\pi(1) \, 1 \, A$ avoids both patterns.

In all subsubcases below let $z=\pi_j$ for some $j >i$.

We start with subsubcase (a).  We must have
$d =\pi_{j+3}$ in order to avoid another occurence
of $(231)$.
Write
$\pi=\pi(1) \, 1 \, A \, z \, x \, y \, d \, \hat{D}$.
We note that for all $\hat{d} \in \hat{D}$ and
all $p \in \pi(1) \, 1 \, A$ we must
have $\hat{d} > z$ and $p<x$ to avoid another $(312)$
or $(231)$ occurence.
Thus, we
have $y=x+1$ and $z=x+2$.  Hence, the value of
$d$ determines the value of $j$ (the position of $z$).
Lastly, we obviously need $\pi(1) \, 1 \, A$ and $\hat{D}$
to be $\{231,312\}$-avoiding.
 By Theorem 2.2,
we now see that we have
$\sum_{j=2}^{n-4}2^{j-2}2^{n-j-4} + 2^{n-5}$ permutations
in this subsubcase.  Hence,  we have
$(n-3)2^{n-6}$ permutations for $n \geq 6$,
one permutation for $n=5$, and none for $n \leq 4$ in
this  subsubcase.

On to subsubcase (b).  We must have $a=\pi_{j-1}$ to
avoid another occurence of $(312)$.
Write
$\pi=\pi(1) \, 1 \, \hat{A} \, a \, z \, x \, y \, D$.
For all $\hat{a} \in \hat{A}$ and for all $d \in D$
we must have $\hat{a} <x$ and $d > z$ in order to
avoid another occurence of either pattern.  Thus, we
have $a=x+1, y=x+2$, and $z=x+3$.  As in
subsubcase (a), we have
$(n-3)2^{n-6}$ permutations for $n \geq 6$, one
permutation for $n=5$, and none for $n \leq 4$ in this
subsubcase.

Next, consider subsubcase (c).
We must have $\pi(1) \, 1 \, A \in S_{j-1}(312;231)$ and
$D \in S_{n-j-2}(\{231,312\})$.
From Theorems 2.2 and 3.4, for each $j \geq 5$ we have
$(j-1)2^{j-6}2^{n-j-3} = (j-1)2^{n-9}$ permutations,
for $j=4$ we have $2^{n-7}$
permutations, and for $j\leq 3$ we have none.
Summing over all valid $j$ we have
$(n-4)(n+1)2^{n-10}$ permutations for $n \geq 7$,
one permutation for $n = 6$, and none
for $n \leq 5$ in this subsubcase.

Lastly, we have subsubcase (d).  A
result similar to that of subsubcase (c) holds. Noting that
$\pi(1) \, 1 \, A \in S_{n-j-2}(\{231,312\})$ and
$D \in S_{j-1}(312;231)$,
from Theorems 2.2 and 3.4, for each $j\leq n-6$ we have
$2^{j-2}(n-j-2)2^{n-j-7} = (n-j-2)2^{n-9}$ permutations.
For $j=n-5$ we have $2^{n-7}$
permutations, and for $j\geq n-4$ we have none.
Summing over all valid $j$ we have
$(n^2-7n+8)2^{n-10}$ permutations for $n \geq 7$,
and none
for $n \leq 6$ in this subsubcase.

Summing over all subsubcases,
we see that we have
$(n^2+19n-70)2^{n-9}$ permutations for $n \geq 6$,
two permutations for $n=5$, and none
for $n \leq 4$
in the subcase $|C|=0$.

Our last case to consider is (iii):
the $(312)$ pattern
straddles 1, i.e. the `3' is in $\pi(1)$, the `2' is
in $\pi(2)$, and 1 serves as the `1' in the pattern.

Let $z1y$ be our $(312)$ pattern and write
$\pi=A\,z\,B\,1\,C\,y\,D$.

We first show that $B$ must be empty.
Assume otherwise and let $b \in B$.
Then we either have another occurence of
$(312)$ with $b1y$ (if $b>y$) or
another occurence of $(312)$ with
$zby$ (if $b<y$).

Next, we show that $|A|+|C| \leq 1$.
Let $a \in A$ and $c \in C$.  We first note
that we must have $c>z$ in order to
avoid another $(312)$ occurence with $z1c$.
We then note that we must have $a<y$ in order
to avoid another $(312)$ occurence with $a1y$.
Hence, for every $a \in A$, $az1$ gives a
$(231)$ occurence, and for every $c \in C$,
$zcy$ gives a $(231)$ occurence.
Thus, $|A|+|C| \leq 1$.

If $|A|=1$ we let $a \in A$ and write
$\pi=a\,z\,1\,y\,D$, where $D \in S_{n-4}(\{231,312\})$.
Further, all elements in $D$ must be larger than $z$
so that we avoid another occurence of $(312)$
with $z$ and $1$.  By Theorem 2.2, we have
$2^{n-5}$ permutations here for $n \geq 5$
and one permutation for $n=4$.

If $|C|=1$ we also have
$2^{n-5}$ permutations for $n \geq 5$
and one permutation for $n=4$
via an argument very similar to that found in the
preceding paragraph.

If $|A|+|C|=0$ we write
$\pi=z\,1\,y\,D$, where
$D \in S_{n-4}(312;231)$ and
again all elements in $D$ are larger
than $z$.  By Theorem 3.4 we have
$(n-3)2^{n-8}$ permutations for $n \geq 7$,
one permutation for $n=6$, and none
for $n \leq 5$ here.

Hence, case (iii) yields
$(n+13)2^{n-8}$ permutations for $n \geq 7$,
five permutations for $n=6$,
two permutations for $n=5$,
and none for $n \geq 4$.

Summing the number of permutations from
all three cases proves the theorem.
\hfill{$\Box$}

For our final class (class (4)) we have the following theorem,
whose proof is more interesting than those above.

{\bf Theorem 5.4}  {\it For $\{\alpha,\beta\} \in
\{\{132,231\},\{132,312\},\{213,231\},\{213,312\}\}$ we
have
\vskip -10pt
\hskip 79pt
$s_n(\emptyset,\{\alpha,\beta\}) = 2^{n-3}$ for
$n \geq 4$.
}

{\it Proof.}   We will use $\{132,312\}$ for our proof.
Let $f_n=s_n(\emptyset;\{132,312\})$.

Let $xzy$ be our $(132)$ pattern and
write $\pi=A\,x\,B\,z\,C\,y\,D$.  First, we show that
$B$ must be empty.  Assume otherwise and let $b \in B$.
Then either $bzy$ or $xby$ is another occurence of
$(132)$ (depending on whether $b<y$ or $b>y$).

Now let $a \in A$.  We must have $y<a<z$ for otherwise
we would have another $(132)$ occurence with $azy$
or more than one $(312)$ occurence with
$axz$ and $axy$.  Also, for $c \in C$ we must have
$c<y$ or else we would have another occurence of
$(132)$ with $xcy$.

We now turn our attention to $D$.  For any $d \in D$ we
must have $d<x$ or $d>z$ in order to avoid
another $(132)$ occurence with $x$ and $z$.  Furthermore,
those elements in $D$ which are larger than $z$ must
be in increasing order so that we avoid another $(132)$
occurence with $x$, and those elements in $D$ which are smaller
than $x$ must be in decreasing order so that we avoid more
than one occurence of $(312)$ with $x$, $y$, and $z$.

Turning back to $A$ and $C$ we now argue that $|A|+|C|=1$.
To see this, note that for any $a \in A$ and any $c \in C$
both $axy$ and $zcy$ are $(312)$ patterns.  Since we
may only have one such pattern we see that $|A|+|C| \leq 1$.
Now assume that both $A$ and $B$ are empty.
With the restrictions on $D$ in the previous paragraph we
see that the pattern $(312)$ is avoided with this assumption.
Hence, $|A|+|B| \geq 1$.

Before putting this all together we note that the above
restrictions show that we have
$\pi=a\,x\,z\,y\,D$ or $\pi=x\,z\,c\,y\,D$ with all elements in $D$
either
smaller than $x$ or larger than $z$.  Hence, the elements preceding D
must
be four consecutive integers which contain both the
patterns $(132)$ and $(312)$ exactly once.

Thus, we have $f_n=f_4 \sum_{i=1}^{n-3} {{n-4} \choose {i-1}} = 2^{n-3}$
permutations in this case (for $n \geq 4$).  This holds since there are
$f_4$ ways to arrange the first four consecutive elements,
we may choose $i=1,2,\dots,n-3$ for the value of
$\min(a,x,y,z)$, and since we are choosing $i-1$
spaces from the $n-4$ spaces
after $y$ in which to place
the decreasing elements of $D$.
\hfill{$\Box$}

{\it Remark}.  We again see another interesting
``mixed restriction" result with
$s_n(132;231) = s_n(\emptyset;\{132,231\})$; i.e.
$(132;231)$ and $(\emptyset;\{132,231\})$ are almost-Wilf equivalent.

\subsection*{\normalsize 5.1.   Generating $S_n(312;123)$ and
$S_n(\emptyset;\{123,312\})$:  On the Almost-Wilf
\vskip -5pt
\hskip 28pt
Equivalence of $(312;123)$ and $(\emptyset;\{123,312\})$}

In this short section we show that the two sets considered are
generated by almost the same rule and let the reader
infer a bijection from these rules.  In the following,
let $n \geq 5$.

Recall (from Section 4.1) that we have defined
$\phi:S_{m-1} \rightarrow S_{m}$ by $\phi(\pi_1 \pi_2 \dots \pi_{m-1})
= (\pi_1+1) (\pi_2+1) \dots (\pi_{m-1}+1) 1$.  We have
also seen the following rule for generating 
$S_n(312;123)$.

{\tt Generating Rule for $S_{n} (312;123)$}:\\
By Theorem 4.6, it is trivial to check
that $S_{n} (312;123) = \{\phi(\pi): \pi \in
S_{n-1}(312;123)\}
\cup \{1 (n-1) n (n-2) (n-3) \dots 3 2, (n-2) (n-1) (n-3) (n-4) \dots 2 1 n\}$.

We now note that we can generate  $S_n(\emptyset;\{123,312\})$ by
the following similar rule.

{\tt Generating Rule for $S_{n} (\emptyset;\{123,312\})$}:\\
By Theorem 5.2, it is trivial to check
that $S_{n} (\emptyset;\{123,312\}) = \{\phi(\pi): \pi \in
S_{n-1}(\emptyset;\{123,$\\
$312\})\}
\cup \{1 n (n-2) (n-1) (n-3) (n-4) \dots 3 2,  (n-1) (n-3) (n-2) \dots 2 1 n \}$.

\subsection*{\normalsize 5.2. Generating $S_n(132;312)$ and
$S_n(\emptyset;\{132,312\})$:  On the Almost-Wilf
\vskip -5pt
\hskip 28pt
Equivalence of $(132;312)$ and $(\emptyset;\{132,312\})$}

In this short section we show that the two sets considered are
generated by exactly the same rule and let the reader
infer a bijection from these rules.  In the following,
let $n \geq 4$.  

From above we have
$\phi:S_{m-1} \rightarrow S_{m}$ by $\phi(\pi_1 \pi_2 \dots \pi_{m-1})
= (\pi_1+1) (\pi_2+1) \dots (\pi_{m-1}+1) 1$.  We also define
$\Phi:S_{m-1} \rightarrow S_{m}$ by $\phi(\pi_1 \pi_2 \dots \pi_{m-1})
= \pi_1 \pi_2 \dots \pi_{m-1} m$.

It is easy to check that the following generation rule
generates both $S_n(132;312)$ and $S_n(\emptyset;\{132,312\})$.
The difference in the sets comes from the initial
sets: $S_4(132;312)=\{3124,4231\}$ and
$S_4(\emptyset;\{132,312\})=\{2413,3142\}$.

{\tt Generating Rule for both $S_n(132;312)$ and $S_n(\emptyset;\{132,312\})$}:\\
To obtain $S_n(\bullet;\bullet)$ from $S_{n-1}(\bullet;\bullet)$ take
 $S_{n}(\bullet;\bullet) = \{\phi(\pi),\Phi(\pi): \pi \in
S_{n-1}(\bullet;\bullet)\}$.

\section*{\normalsize 6. Summary and Questions}

Below we give a table summarizing the above
results and present some remaining questions.
The top half of the table's results comes from
Section 3, and the bottom half comes from
Section 4.

\begin{center}

\begin{tabular}{||l|l||} \hline\hline
\multicolumn{1}{||c|}{Almost-Wilf Class, $\cal{W}$}   &
\multicolumn{1}{c||}{$s_n(T), T \in \cal{W}$}\\ \hline
   &\\
$A = \overline{(123;321)}$&0  for $ n \geq 6$\\
&\\
$B=\overline{(123;132)}$&$(n-2)2^{n-3}$ for $n \geq 3$\\
&\\
$C=\overline{(123;231)}$&$2n-5$ for $n \geq 3$\\
&\\
$D=\overline{(132;213)}$&$n2^{n-5}$ for $n \geq 4$\\
&\\
$E=\overline{(132;231)}$&$2^{n-3}$ for $n \geq 3$\\
&\\
\hline
&\\
$F=\overline{(\emptyset;\{123,\!321\})}$&$0$ for $n \geq 6$\\
&\\
$G=\overline{(\emptyset;\{123,\!231\})}$&$2n-5$ for $n \geq 5$\\
&\\
$H=\overline{(\emptyset;\{123,\!132\})}$&${{n-3} \choose 2}2^{n-4}$ for
$n
\geq 5$ \\
&\\
$I=\overline{(\emptyset;\{132,\!213\})}$&$(n^2+21n-28)2^{n-9}$ for $n
\geq
7$\\
&\\
$J=\overline{(\emptyset;\{132,\!231\})}$&$2^{n-3}$ for $n \geq 4$\\
&\\
\hline\hline
\end{tabular}

\end{center}
\vskip 10pt
{\small

\centerline{ The Cooresponding Almost-Wilf Classes}

A.  $\{(123;\!321), (321;\!123)  \}$

B.  $\{(123;\!132), (123;\!213), (132;\!123), (213;\!123),
(231;\!321), (312;\!321), (321;\!231),(321;\!312)\}$

C.  $\{(123;\!231),(123;\!312),(132;\!321),(213;\!321),(231;\!123),
(312;\!123),(321;\!132),(321;\!213)\}$

D.  $\{(132;\!213), (213;\!132), (231;\!312), (312;\!231)  \}$

E. $\{(132;\!231),(132;\!312), (213;\!231),(213;\!312), (231;\!132),
(231;\!213), (312;\!132),(312;\!213)\}$

F.  $\{(\emptyset;\{123,\!321\}), (\emptyset;\{321,\!123\})  \}$

G.  $\{(\emptyset;\{123,\!231\}),(\emptyset;\{123,\!312\}),
(\emptyset;\{132,\!321\}),(\emptyset;\{213,\!321\})\}$

H.  $\{(\emptyset;\{123,\!132\}),(\emptyset;\{123,\!213\}),
(\emptyset;\{231,\!321\}),(\emptyset;\{312,\!321\})\}$

I.  $\{(\emptyset;\{132,\!213\}),(\emptyset;\{231,\!312\})\}$

J.  $\{(\emptyset;\{132,\!231\}),(\emptyset;\{132,\!312\}),
(\emptyset;\{213,\!231\}),(\emptyset;\{213,\!312\})\}$
}

We would like very much to see formulas for
$s_n(\emptyset;\{(123)^2\})$ and $s_n(\emptyset;\{(132)^2\})$
determined to
finish the study of
$s_n(\emptyset;\{\alpha,\beta\})$ for all
$\alpha,\beta \in S_3$.
Note that B\'ona, in [B2], has given the generating function and
a recursive formula for the sequence
$\{s_n(\emptyset;\{(132)^2\})\}_n$, however a formula for
$s_n(\emptyset;\{(132)^2\})$ is not immediate.

{\bf Acknowledgment}

I would like to thank Alek Vainshtein for calling my
attention to some relevant references.

{\bf References}
\footnotesize
\parskip=5pt

[B1] M. B\'ona, Exact Enumeration of $1342$-avoiding
Permutations; A Close Link with Labelled Trees and
Planar Maps, {\it Journal of Combinatorial
Theorey, Series A}, {\bf 175} (1997), 55-67.

[B2] M. B\'ona, Permutations Avoiding Certain Patterns.  The
Case of Length 4 and Some Generalizations,
{\it Discrete Mathematics}, {\bf 80} (1997), 257-272.

[B3] M. B\'ona, The Permutation Classes Equinumerous to the
Smooth Class, {\it Electronic Journal of Combinatorics},
{\bf 5(1)} (1998), R31.

[B4] M. B\'ona, Permutations with One or Two
$132$-sequences, {\it Discrete Mathematics}, {\bf 181} (1998), 267-274.

[BLPP1] E. Barcucci, A. Del Lungo, E. Pergola, and
R. Pinzani, From Motzkin to Catalan Permutations,
{\it Discrete Mathematics}, {\bf 217} (2000), 33-49.

[BLPP2] E. Barcucci, A. Del Lungo, E. Pergola, and
R. Pinzani, Permutations Avoiding an Increasing Number of
Length-Increasing Forbidden Subsequences, {\it Discrete
Mathematics and Theoretical Computer Science},
{\bf 4} (2000), 31-44.

[C] E. Catalan, Note sur une \'Equation aux Differences Finies,
{\it Journal des Mathematiques Pures et Appliqu\'es},
{\bf 3} (1838), 508-516.

[CW] T. Chow and J. West, Forbidden Subsequences and Chebyshev
Polynomials, {\it Discrete Mathematics}, {\bf 204} (1999), 119-128.

[Ge] I. Gessel, Symmetric Functions and P-recursiveness,
{\it Journal of Combinatorial Theory, Series A}, {\bf 53}
(1990), 257-285.

[Gu] O. Guibert, Combinatoires des Permutations a Motifs Exclus en
Liaison
avec Mots,
Cartes Planaires et Tableaux de young, {\it Th\`ese de
l'Universit\'e de Bordeaux I} (1995).

[JR] M. Jani and R. Rieper, Continued Fractions and Catalan Problems,
{\it Electronic Journal of Combinatorics}, {\bf 7(1)} (2000), R45.

[Kn] D. Knuth, \underline{The Art of Computer Programming}, vol. 3,
Addison-Wesley, Reading, MA, 1973.

[Kr] D. Kremer, Permutations with Forbidden Subsequences and a
Generalized Schr\"oder Number, {\it Discrete Mathematics},
{\bf 218} (2000), 121-130.

[Krt] C. Krattenthaler, Restricted Permutations, Continued Fractions,
and Chebyshev Polynomials, (preprint 2000) math.CO/0002200.

[M] T. Mansour, Permutations Avoiding a Pattern from $S_k$
and at Least Two Patterns from $S_3$, (preprint 2000)
math.CO/0007194.

[MV1] T. Mansour and A. Vainshtein, Restricted Permutations, Continued
Fractions, and Chebyshev Polynomials, {\it Electronic
Journal of Combinatorics}, {\bf 7(1)} (2000), R17.

[MV2] T. Mansour and A. Vainshtein, Restricted $132$-Avoiding
Permutations, (preprint 2000)\\
math.CO/0010047.

[MV3] T. Mansour and A. Vainshtein, Restricted Permutations and
Chebyshev Polynomials, (preprint 2000)
math.CO/0011127.

[N] J. Noonan, The Number of Permutations Containing Exactly
One Increasing Subsequence of Length Three, {\it Discrete Mathematics},
{\bf 152} (1996), 307-313.

[NZ] J. Noonan and D. Zeilberger, The Enumeration of Permutations with a
Prescribed Number of ``Forbidden" Patterns, {\it Advances in Applied 
Mathematics}, {\bf 17} (1996), 381-407.

[R] A. Robertson, Permutations Containing and Avoiding
$123$ and $132$ Patterns, {\it Discrete Mathematics
and Theoretical Computer Science}, {\bf 4} (1999),
151-154.

[Ri] D. Richards, Ballot sequences and Restricted Permutations,
{\it Ars Combinatoria}, {\bf 25} (1988), 83-86.

[RWZ] A. Robertson, H. Wilf, and D. Zeilberger,
Permutation Patterns and Continued Fractions,
{\it Electronic Journal of Combinatorics},
{\bf 6(1)} (1999), R38.

[S] Z. Stankova, Forbidden Subsequences, {\it Discrete Mathematics},
{\bf 132} (1994), 291-316.

[SS] R. Simion and F. Schmidt, Restricted Permutations,
{\it European Journal of Combinatorics} {\bf 6} (1985), 383-406.

[W1] J. West, Ph.D. Thesis.

[W2] J. West, Sorting Twice Through a Stack, {\it Theoretical
Computer Science}, {\bf 117} (1993), 303-313.

\end{document}